\documentclass[12pt]{amsart}
 \usepackage{tikz}
\usepackage{amssymb,amsmath}
\def\Xint#1{\mathchoice
{\XXint\displaystyle\textstyle{#1}}%
{\XXint\textstyle\scriptstyle{#1}}%
{\XXint\scriptstyle\scriptscriptstyle{#1}}%
{\XXint\scriptscriptstyle%
\scriptscriptstyle{#1}}%
\!\int}
\def\XXint#1#2#3{{\setbox0=\hbox{$#1{#2#3}{%
\int}$ }
\vcenter{\hbox{$#2#3$ }}\kern-.6\wd0}}

\def\dashint{\Xint-}

\newtheorem{theorem}{Theorem} 

\numberwithin{equation}{section}

\def\di{\partial}
\def\dib{\bar\partial}
\numberwithin{equation}{section}
\def\simleq{\underset\sim<}

\def\T{\text}
\def\1#1{\overline{#1}}
\def\2#1{\widetilde{#1}}
\def\3#1{\widehat{#1}}
\def\4#1{\mathbb{#1}}
\def\5#1{\frak{#1}}
\def\6#1{{\mathcal{#1}}}
\def\C{{\4C}}

\def\B{\Bbb B}

\def\phi{\varphi}

\newtheorem{Thm}{Theorem}
\newtheorem{Cor}[Thm]{Corollary}
\newtheorem{Pro}[Thm]{Proposition}
\newtheorem{Lem}[Thm]{Lemma}
\theoremstyle{definition}\newtheorem{Def}[Thm]{Definition}
\theoremstyle{remark}
\newtheorem{Rem}[Thm]{Remark}
\newtheorem{Exa}[Thm]{Example}

\def\Label#1{\label{#1}}
\def\bl{\begin{Lem}}
\def\el{\end{Lem}}
\def\bp{\begin{Pro}}
\def\ep{\end{Pro}}
\def\bt{\begin{Thm}}
\def\et{\end{Thm}}
\def\bc{\begin{Cor}}
\def\ec{\end{Cor}}
\def\bd{\begin{Def}}
\def\ed{\end{Def}}
\def\br{\begin{Rem}}
\def\er{\end{Rem}}
\def\be{\begin{Exa}}
\def\ee{\end{Exa}}
\def\bpf{\begin{proof}}
\def\epf{\end{proof}}
\def\ben{\begin{enumerate}}
\def\een{\end{enumerate}}

\def\1alpha{[\frac1\alpha]}
\def\T{\text}

\def\C{{\Bbb C}}

\numberwithin{equation}{section}
\def\T{\text}
\newcommand{\om}{\omega}

\newcommand{\no}[1]{\|{#1}\|}

\newtheorem{definition}[theorem]{Definition }
\newtheorem{lemma}[theorem]{Lemma  }
\newtheorem{proposition}[theorem]{Proposition  }
\newtheorem{corollary}[theorem]{Corollary }
\newtheorem{example}[theorem]{\it Example }

\begin{document}
\title[Extension of $L^2$, $\dib$-closed, forms]
{Extension of $L^2$, $\dib$-closed, forms}
\author[ L.~Baracco, S.~Pinton and G.~Zampieri]
{Luca Baracco, Stefano Pinton and Giuseppe Zampieri}
\address{Dipartimento di Matematica, Universit\`a di Padova, via 
Trieste 63, 35121 Padova, Italy}
\email{baracco@math.unipd.it, pinton@math.unipd.it, zampieri@math.unipd.it}
\maketitle
\begin{abstract}
We prove extension of a $\dib$-closed, smooth,  form from the intersection of a pseudoconvex domain with a complex hyperplane to the whole domain. The extension
form is $\dib$-closed, has harmonic coefficients and its $L^2$-norm is estimated by the $L^2$-norm of the trace. For holomorphic functions this is proved by Ohsawa-Takegoshi \cite{OT87}. For forms of higher degree, this is stated in \cite{M93}. It seems, however, that the proof contains a gap because of the use of a a singular  weight and the failure of regularity for the solution of the related $\dib$-equation. 

There is a rich literature on the subject (cf. among otheres \cite{D97}, \cite{S96}) but it does not  seem to contain complete answer to the question. Also, the problem of extending cohomology classes of $\dib$ of higher degree in a compact K\"ahler space is addressed in \cite{K11} and \cite{B11}. Apart from the formal analogy, this has little in common with our problem in which these classes are $0$.  We also believe, comparing to the preceding literature, that our approach is original and, somewhat, simpler.  

\end{abstract}

We recall 
\bt
\Label{t1.1}
Let $D\subset\subset\C^n$ be a bounded smooth pseudoconvex domain of diameter $\le1$, $\alpha$ a $\dib$-closed form of degree $\ge 1$ such that $\alpha_J=0$ for $1\notin J$
and supp$\,\alpha\subset\{z\,:\,|z_1<\delta\}$.
 Then there is a solution $u=u_\delta$ in $L^2$ to the problem
\begin{equation}
\Label{1.1}
\begin{cases}
\dib u=\alpha,
\\
\no{z_1 u}_0\le c\delta\no{\alpha}_0,\quad\T{for $c$ independent of $\delta$, $\alpha$ and $D$}.
\end{cases}
\end{equation}
\et
We refer to \cite{BPZ15} for a proof which only relies on the  Kohn-H\"ormander-Morrey estimates in the weighted $L^2$ space; moreover, ``selfboundedness" of the gradient of the weights is never used. 
The problem is to extend a $\dib$-closed, smooth, form $f$ from the slice $D^0=D\cap \{z:\,z_1=0\}$ to the full $D$ with control of the $L^2$ norm.  With the notation  $\4{D}_\delta$ for  the $\delta$-disc in the $z_1$-plane, this can be achieved  by taking a pseudoconvex approximation $D_\nu\nearrow D$, choosing $\delta=\delta_\nu$ such that the $\delta$-strip 
$\4{D}_\delta\times D^0_{\nu+1}$ contains
 $D_\nu\cap(\4{D}_\delta\times\C^{n-1})$
taking a family of cut-offs $\chi_\delta(z_1)$ with unitary mass which is $1$ for $|z_1|\le\delta$ and $0$ for $|z_1|>2\delta$, and a family of solutions $\{u_\delta\}$ of \eqref{1.1} on $D_\nu$  for the choice of the form $\alpha_\delta=\frac{\dib\chi_\delta(z_1) f(z_2,...,z_n)}{z_1}$. Each form $f_\delta:=-z_1u_\delta+\chi_\delta f$ is $\dib$-closed in $D_\nu$
and the family $\{\no{f_\delta}_0\}$ is uniformly bounded by $c\no{f}_0^{D^0}$. Then there is a subsequential $L^2$-weak limit $\tilde f$ on $D$, which satisfies $\dib \tilde f=0$ and 
\begin{equation}
\Label{1.1,5}
\no{\tilde f}^D_0\simleq \no{f}_0^{D^0}.
\end{equation}
 If the degree of the forms is $0$, that is, if $f_\delta$ and $\tilde f$ are holomorphic functions, then
\begin{equation}
\Label{1.2}
f_\delta\to \tilde f\T{ pointwise}.
\end{equation}
This can be readily checked recalling that, over holomorphic functions, weak convergence implies pointwise convergence. On the other hand, $\dib$ is elliptic on functions and thus $u_\delta\in C^\infty$, $z_1u_\delta|_{z_1=0}=0$ and therefore $f_\delta|_{z_1=0}=f$; it follows
\begin{equation}
\Label{1.3}
\tilde f|_{z_1=0}=f.
\end{equation}
 Thus holomorphic extension with the estimate \eqref{1.1,5} is proved, for functions. What follows is dedicated to first show that \eqref{1.2} remains true, for an accurate choice of the $u_\delta$'s, in general degree. 
And next to prove a slightly weaker version of \eqref{1.3}, that is,
\begin{equation}
\Label{1.3bis}
j^*\tilde f=f\quad\T{where $j:\,D^0\hookrightarrow D$ is the embedding}.
\end{equation}
For the first, we have to take a minimizer with respect to the norm $\int_{D_\nu}|z_1|^2|\cdot|^2dV$
of the affine space of solutions of \eqref{1.1}. The minimizer
is a $L^2$ limit $z_1u^\mu\to z_1 u_\delta$; in particular, it inherits the estimate, uniform in $\delta$, $\no{z_1u_\delta}_0^D\simleq  \no{f}^{D^0}_0$. 
Moreover,
 $u_\delta$ does not necessarily solve $\dib u_\delta=\alpha_\delta$ but it certainly solves $\dib (z_1u_\delta)=z_1\alpha_\delta$. Also, since it is a minimizer, it satisfies, for every $v\in C^\infty_c$
\begin{equation*}
\begin{split}
0&= (|z_1|^2u_\delta,\dib v)_0
\\
&=(z_1u_\delta,\dib(z_1v))_0\\
&=(\vartheta(z_1u_\delta),z_1v).
\end{split}
\end{equation*}
Thus $z_1u\in\ker\vartheta$ for $z_1\neq0$. But, in fact, also for $z_1=0$. In fact
for $\psi\in C^\infty_c$, we decompose $\psi=\psi^++\psi^\epsilon$ where $\psi^+ $ is supported by $|z_1|>\epsilon$ and $\psi^\epsilon$ by $|z_1|<2\epsilon$. Then
\begin{equation*}
\begin{split}
(\vartheta(z_1u), \psi)&=(\vartheta(z_1u),\psi^++\psi^\epsilon)_0
\\
&=(\vartheta(z_1u), \psi^\epsilon)_0
\\
&=(z_1u,\dib \psi^\epsilon)_0
\\
&\underset{\T{Cauchy-Schwarz}}\le\no{\chi_{2\epsilon}z_1u}_0\no{\dib \psi^\epsilon}_0
\\
&\to 0,
\end{split}
\end{equation*}
since $z_1u\in L^2$ and $\no{\dib \psi^\epsilon}_0$ is uniformly bounded. 
Thus $\vartheta(z_1u_\delta)=0$ and in particular $z_1u_\delta\in C^\infty$.
Now, we remark that $\tilde f$ is the limit not only of $-z_1u_\delta+\chi_\delta f$ but also of $-z_1u_\delta$. Thus, not only it satisfies $\dib \tilde f=0$, but also $\vartheta \tilde f=0$. It follows
$$
\Delta \tilde f=0;
$$
in particular, $\tilde f\in C^\infty$.
As for $f_\delta$, we have $\dib f_\delta=0$ but,  we only have $\vartheta f_\delta=\vartheta\chi_\delta f$, which yields $\vartheta f_\delta\to 0$ ($L^2$-weakly) and $f_\delta\in C^\infty$.
To carry on our proof, a more subtle analysis is needed. We remark that
$$
\vartheta \dib f_\delta=0,
$$
and
\begin{equation*}
\begin{split}
\dib\vartheta f_\delta&=\dib\vartheta(\chi_\delta f)
\\
&\underset{\T{(*)}}=\dib(\chi_\delta \vartheta f)
\\
&=\sum_i\dib_i\chi_\delta\sum_j\di_jf_{jK}+\chi_\delta\sum_i(\sum_j\dib_i\di_jf_{jK}),
\end{split}
\end{equation*}
where (*) is explained by the fact that $\chi_\delta$ depends on $z_1$ and $f$ on $z_2,...,z_n$, only.
It follows
\begin{equation}
\Label{1.5}
\begin{split}
|\Delta f_\delta|&\le c|\dib\chi_\delta|+c\chi_\delta
\\
&\le c\frac{\chi_{2\delta}}{\delta}.
\end{split}
\end{equation}
In particular, $f_\delta$ is harmonic for $|z_1|>2\delta$ and $f_\delta\to \tilde f$ uniformly on compact subsets of $z_1\neq0$. 
Thus \eqref{1.2} holds for $z_1\neq0$. As for $z_1=0$, with the notation $\sigma_{2n}$ for the volume of the unit ball in $\C^n$, we have
\begin{equation}
\Label{1.50}
\begin{split}
\dashint_{\B_r(z_o)}f_\delta dV-f_\delta(z_o)&=\frac1{\sigma_{2n}r^{2n}}\int_0^r\int_{b\B_t(0)}[f_\delta(z_o+\xi)-f_\delta(z_o)]dS_\xi dt
\\
&=\frac1{\sigma_{2n}r^{2n}}\int_0^r t\int_{b\B_t(0)}\int_0^1\nabla f_\delta(z_o+s\xi)\cdot \frac\xi t ds dS_\xi dt
\\
&
\underset{\T{Divergence Th.}}=\frac1{\sigma_{2n}r^{2n}}\int_0^rt\int_0^1s\int_{\B_t(0)}\Delta f_\delta(z_o+s\xi)dV_\xi ds dt
\\
&\underset{\T{$s\xi=\zeta$}}=\frac1{\sigma_{2n}r^{2n}}\int_0^r t\int_0^1 s^{-2n+1}\int_{\B_{st}(z_o)}\Delta f_\delta dV_\zeta ds dt.
\end{split}
\end{equation}
Recall that $|\Delta f_\delta|\le c\frac{\chi_\delta}\delta$ and hence, for $z_o$ in the plane $z_1=0$,
\begin{equation}
\Label{1.51}
\begin{split}
\int_{\B_{st}(z_o)}|\Delta f_\delta|dV&\simleq\frac1\delta \T{Vol}((\4{D}_\delta\times \B^{2n-2}_{st})\cap \B^{2n}_{st})
\\
&\le (st)^{2n-1}.
\end{split}
\end{equation}
Combination of \eqref{1.50} and \eqref{1.51} yields, for $z_o$ in the plane $z_1=0$,
\begin{equation}
\Label{1.6}
\begin{split}
\Big|\dashint_{\B_r(z_o)}f_\delta dV-f_\delta(z_o)\Big |&\simleq r^{-2n}\int_0^r t^{2n} dt \int_0^1ds
\\
&\simleq r.
\end{split}
\end{equation}
We know that we have $L^2$ weak convergence $f_\delta\to \tilde f$; in particular, by the choice of the characteristic function of $\B_r(z_o)$ as the test function, we have
\begin{equation}
\Label{1.7}
\lim_\delta \Big|\dashint_{\B_r(z_o)}f_\delta dV-\dashint_{\B_r(z_o)}\tilde f dV\Big|=0,\quad \T{$r$ fixed}.
\end{equation}
Finally, by the harmonicity of $\tilde f$
\begin{equation}
\Label{1.8}
\dashint_{\B_r(z_o)}\tilde f dV-\tilde f(z_o)\equiv 0\quad\T{ for any $r$.}
\end{equation}
Plugging together \eqref{1.6}, \eqref{1.7} and \eqref{1.8}, we get, for $z_o$ in the plane $z_1=0$
\begin{equation}
\Label{1.9}
\lim_\delta (f_\delta(z_o)-\tilde f(z_o))=0.
\end{equation}
(In fact, for any $\epsilon$, we first choose $r$ such that  \eqref{1.6} is $<\epsilon$; under this choice of $r$, we next choose $\delta$ such that \eqref{1.7} is also $<\epsilon$.)
Thus we have proved \eqref{1.2} for a form of general degree; we pass now to \eqref{1.3bis}. This is obtained as an immediate consequence of $j^*(z_1u_\delta)=0$, that is
\begin{equation}
\Label{1.10}
(z_1u_\delta)_H|_{z_1=0}=0\quad\T{ for any $H$ which does not contain $1$}.
\end{equation}
We point out that this is by no means evident because we do not know whether $u_\delta$ is smooth or $ L^2$. We also point out that our method does not yield $z_1u_\delta|_{z_1=0}=0$ in full, but just for the ``tangential" part (the one which collects multiindices which do not contain 1). Thus our program is to prove \eqref{1.10}. We decompose
\begin{equation}
\Label{1.11}
u_\delta=\frac h{z_1}+g,\quad \T{ for $h,\,g\,\in C^\infty$, $h+z_1g\in L^2$.}
\end{equation}
We reason by contradiction and assume that for an index $J$ with $ 1\notin J$, that we also write $J=iK$, we have $h_{iK}\neq0$.
We choose a family of test functions
\begin{equation}
\Label{1.12}
\psi_\epsilon=\chi(\frac{|z_1-3\epsilon|^2}{\epsilon^2})\chi(|z_2|)...\chi(|z_n|);
\end{equation}
we also use the notation $\psi(z'):=\chi(|z_2|)...\chi(|z_n|)$. We arrange coordinates so that supp$\psi_\epsilon\subset\subset D$; we notice that supp$\psi_\epsilon\subset\subset\{z:\,z_1\neq0\}$. Let $z_1u^\mu\to z_1u_\delta$ be the approximation of the minimizer $u_\delta$ in $\int|z_1|^2|\cdot|^2\,dV$ norm. We have, for any $\psi\in Dom(\dib^*)$ such that $\dib^*\psi|_{z_1=0}=0$,
\begin{equation}
\Label{1.13}
\begin{split}
|(\dib(u^\mu-u_\delta),\psi)|&\le \left(\int |z_1(u^\mu-u_\delta|^2\,dV\right)^{\frac12}\left(\int\Big|\frac{\dib^*\psi}{z_1}\Big|^2\,dV\right)^{\frac12}
\\
&\to0.
\end{split}
\end{equation}
 Remember also that $\dib u^\mu\equiv \frac{\dib \chi_\delta f}{z_1}\in C^\infty$ (independent of $\mu$) that we also call $\alpha$. Thus, by \eqref{1.13} 
 applied for $\psi=\psi_\epsilon \bar \om_H$, with $H=1 iK$, we get
\begin{equation}
\Label{1.14}
\lim_\epsilon\int(\dib u_\delta)_H\psi_\epsilon\,dV=\lim_\epsilon\int\alpha_H\psi_\epsilon\,dV=0.
\end{equation}
Now,
\begin{equation}
\begin{split}
\Label{1.15}
\int&\di_{\bar z_1}(u_\delta)_{iK}\psi_\epsilon=-\int(u_\delta)_{iK}\frac{\bar z_1-3\epsilon}{\epsilon^2}\dot \chi(\frac{|z_1-3\epsilon|^2}{\epsilon^2})\psi(z')dV
\\
&\underset{\T{\eqref{1.11}}}=-\int h_{iK}(z')\psi(z')\,dV_{z'}\Big(\frac1{\epsilon^2}\int_{\Delta_{2\epsilon}(3\epsilon)}\dot \chi(\frac{|z_1-3\epsilon|^2}{\epsilon^2})\frac{\bar z_1-3\epsilon}{z_1}dx_1dy_1\Big)
\\
&\hskip2cm-\int g_{iK}(0,z')\psi(z')\,dV_{z'}\Big(\frac1{\epsilon^2}\int_{\Delta_{2\epsilon}(3\epsilon)}\dot \chi(\frac{|z_1-3\epsilon|^2}{\epsilon^2})(\bar z_1-3\epsilon)dx_1dy_1\Big)+O(\epsilon^2)
\\
&=-\int h_{iK}(z')\psi(z')\,dV_{z'}\Big(\frac1{\epsilon^2}\int_{\Delta_{2\epsilon}(3\epsilon)}\dot \chi(\frac{|z_1-3\epsilon|^2}{\epsilon^2})\frac{\bar z_1-3\epsilon}{z_1}dx_1dy_1\Big)+O(\epsilon)+O(\epsilon^2).
\end{split}
\end{equation}
We show that the integral in the last line is a constant different from $0$. Under the coordinate change  in $\C$, $z_1-3\epsilon=\rho\zeta,\,\,\rho\in[0,2\epsilon],\,\zeta\in\di\Delta$, we have
\begin{equation}
\begin{split}
\int_{\Delta_{2\epsilon}(3\epsilon)}\dot \chi(\frac{|z_1-3\epsilon|^2}{\epsilon^2})\frac{\bar z_1-3\epsilon}{z_1}dx_1dy_1&=\int_0^{2\epsilon}\dot\chi(\frac{\rho^2}{\epsilon^2})\int_{\di\Delta}\frac{\rho\bar \zeta}{\rho\zeta+3\epsilon}\frac\rho{i\zeta} d\zeta d\rho
\\
&=\int_0^{2\epsilon}\dot\chi(\frac{\rho^2}{\epsilon^2})\frac{\rho^2}{3i\epsilon}\int_{\di\Delta}\frac1{\zeta^2(1+\frac{\rho\zeta}{3\epsilon})}d\zeta d\rho
\\
&\underset{\T{Residue Th.}}=\int_0^{2\epsilon}\dot\chi(\frac{\rho^2}{\epsilon^2})\frac{\rho^2}{3i\epsilon}2\pi i\left(\frac{-\rho}{3\epsilon}\right)d\rho
\\
&=-\frac{2\pi}{9\epsilon^2}\int_0^{2\epsilon}\dot\chi(\frac{\rho^2}{\epsilon^2})\rho^3d\rho
\\
&\underset{\T{by the change $t=\frac{\rho^2}{\epsilon^2}$}}=-\frac{2\pi}{9\epsilon^2}\int_0^4\dot\chi(t)t\frac{\epsilon^4}2 dt
\\
&=-\frac{\pi\epsilon^2}9[\Big|_0^4 \chi(t)t-\int_0^4\chi(t)dt]
\\
&=\frac{\pi\epsilon^2}9.
\end{split}
\end{equation}
Thus, if $h_{iK}\neq0$, then the limit for $\epsilon\to0$ of \eqref{1.15} is not $0$. On the other hand, whatever the value of $h_{1K}$ is, we have, when $i\neq1$,
\begin{equation}
\Label{1.16}
\begin{split}
\int\di_{\bar z_i}(u_\delta)_{1K}\psi_\epsilon dV&=-\int (u_\delta)_{1K}\chi\left(\frac{|z_1-3\epsilon|^2}{\epsilon^2}\right)\di_{\bar z_i}\chi(z') dV
\\
&=-\int h_{1K}(z')\di_{\bar z_i}\chi(z')\int_{\Delta_{2\epsilon}(3\epsilon)}\frac{\chi\left(\frac{|z_1-3\epsilon|^2}{\epsilon^2}\right)}{z_1}dx_1dy_1dV_{z'}+O(\epsilon^2)
\\
&=-\int h_{1K}(z')\dot\chi(z')O(\epsilon)+O(\epsilon^2)
\\
&=O(\epsilon).
\end{split}
\end{equation}
(Clearly, in case $(u_\delta)_{1K}$ is bounded, that is $h_{1K}=0$, we have that \eqref{1.16} is indeed $O(\epsilon^2)$.) In any case, for $ iK$ not containing $1$ and with $H=1iK$, from the identity
$(\dib u_\delta)_H=\di_{\bar z_1}((u_\delta)_{iK})-\di_{\bar z_i}((u_\delta)_{1K})$, we get, combining \eqref{1.15} and \eqref{1.16},
$$
\int(\dib u_\delta)_{1iK}\psi_\epsilon dV=c+0(\epsilon),\quad c\neq0,
$$
which contradicts \eqref{1.14}.
Summarizing up, we have proved that in the decomposition $(u_\delta)_{iK}=\frac {h_{iK}}{z_1}+g_{iK} $, we cannot have $h_{iK}\neq0$ when $1\notin iK$. Thus, $(u_\delta)_{iK}$ is smooth and in particular $z_1(u_\delta)_{iK}|_{z_1=0}=0$, and hence $(f_\delta)_{iK}|_{z_1=0}=f_{iK}$. If we plug into \eqref{1.9}, we get   $\tilde f_{iK}|_{z_1=0}=f_{iK}$.
All in all, we have got
\bt
\Label{t1.2}
Let $D$ be a smooth pseudoconvex domain of $\C^n$, let $D^0$ be the slice $D\cap \{z\,:\,z_1=0\}$ and let $j\,:\,D^0\to D$ be the embedding.
For every $L^2\cap C^\infty$-form $f$ on $D^0$ such that $\dib f=0$,  there is an $L^2\cap C^\infty$ form $\tilde f$ on $D$ which is $\dib$-closed, has harmonic coefficients and satisfies
\begin{equation*}
\begin{cases}
j^*\tilde f=f,
\\
\no{\tilde f}_0^D\le c\no{f}_0^{D^0}.
\end{cases}
\end{equation*}
\et


\begin{thebibliography}{McNV15}
\bibitem{BPZ15}{\bf L. Baracco}---Extension of $L^2$ holomorphic functions, (2015)
\bibitem{B96}{\bf B. Berndtsson}---The extension theorem of Ohsawa-Takegoshi and the theorem of Donnelly-Fefferman, {\em Ann. Inst. Fourier} {\bf 46} n. 4 (1996), 1083--1094
\bibitem{B11} {\bf B. Berndtsson}---$L^2$ etxension of $\dib$-closed forms, arxiv:1104.4620v1
\bibitem{B13}{\bf Z. Blocki}---Suita conjecture and the Ohsawa-Takegoshi extension theorem, {\em Inv. Math.} {\bf 193} (2013), 149--158
\bibitem{B14}{\bf Z. Blocki}---Cauchy-Riemann meet Monge-Amp\`ere, {\em Bull. Math. Sci.} {\bf 4} (2014), 433--480
\bibitem{C11}{\bf B.Y. Chen}---A simple proof of the Ohsawa-Takegoshi extension theorem, arxiv : 11052430v1
\bibitem{D97} {\bf J.P. Demailly}---On the Ohsawa-Takegoshi-Manivel $L^2$ extension Theorem, {\em  Prog. Math} {\bf 188}, 47--82
\bibitem{K11}{\bf V. Koziarz}---Extension with estimates of cohomology classes, {\em Man. Math.} {\bf 134} n. 1-2 (2011), 43--58




\bibitem{M93} {\bf L. Manivel}---Un th\'eor\`eme de prolongement $L^2$ des sections holomorphes d'un fibr\'e hermitien, {\em Math. Z.} {\bf 212} n. 1 (1993), 107--122
\bibitem{McNV15}{\bf J.D. McNeal and D. Varolin}---$L^2$ extension of $\dib$-closed forms from a hypersurface, arxiv:1502.08054v1
\bibitem{O88} {\bf T. Ohsawa}---On the extension of $L^2$ holomorphic functions II, {\em Publ. RIMS} {\bf 24} (1988), 265--275
\bibitem{OT87}{\bf T. Ohsawa and K. Takegoshi}---On the extension of $L^2$, holomorphic functions, {\em Math. Z.} {\bf 195} (1987), 197--204
\bibitem{P15}{\bf S. Pinton}---$L^2$-estimates with multipliers for holomorphic extension, (2015)
\bibitem{S96} {\bf Y.T. Siu}---The Fujita conjecture and the extension theorem of Ohsawa-Takegoshi, {\em Geometric Complex Analysis, Hayama, World Sci. Publ.} (1996)
\bibitem{S11}{\bf Y.T. Siu}---Section extension from hyperbolic geometry of punctured disk and holomorphic families of flat bundles, arxiv:11.042563
\bibitem{S10} {\bf E. Straube}---Lectures on the $L^2$-Sobolev theory of the $\bar\partial$-Neumann problem, {\em ESI Lect. in Math. and Physics} (2010)
 
\end{thebibliography}
\end{document}